\documentclass[graybox]{svmult}

\usepackage{amsfonts}
\usepackage{amsmath}
\usepackage{a4}
\usepackage{url}
\usepackage{graphicx}

\newcommand{\Z}{\mathbb{Z}}
\newcommand{\R}{\mathbb{R}}
\newcommand{\C}{\mathbb{C}}

\DeclareMathOperator{\GL}{GL}
\DeclareMathOperator{\conv}{conv}

\begin{document}

\title*{Exploiting symmetries in polyhedral computations}

\author{Achill Sch\"urmann}
\institute{Achill Sch\"urmann \at 
Institute of Mathematics, University of Rostock, 18051 Rostock, Germany,
\email{achill.schuermann@uni-rostock.de}}


\maketitle

\abstract{In this note we give a short overview on symmetry exploiting techniques 
in three different branches of polyhedral computations: 
The representation conversion problem,
integer linear programming and lattice point counting.
We describe some of the future challenges and 
sketch some directions of potential developments.
}

\section{Introduction}

Symmetric polyhedra such as the Platonic and Archimedean solids
have not only fascinated mathematicians since time immemorial.
They occur frequently in diverse contexts of art and science.
Less known to a general audience, but of great importance
to modern mathematics and its applications, are higher dimensional
analogues of these familiar objects.
One {\em standard description} is as 
a set of solutions to a system of
linear inequalities 
\begin{equation*} 
P =
\{
x\in \R^n :
Ax \leq b
\}
,
\end{equation*}
where~$A$ is a real ${m\times n}$ matrix and $b\in\R^m$.
A prominent example is the {\em $n$-cube} obtained
by $2n$ inequalities $\pm x_i \leq 1$.
It has $2^n$ {\em vertices} ({\em extreme points}) with coordinates $\pm 1$
and its {\em group of symmetries} is the {\em hyperoctahedral group} of
order $2^n n!$.  

\medskip

Linear models, and therefore polyhedra, are used in a wide range of 
mathematical problems and in applications such as
transportation logistics,
machine scheduling, time tabling, air traffic flow management
and portfolio planning.
They are central objects 
in Mathematical Optimization (Mathematical Programming)
and are for instance heavily used in Combinatorial Optimization.
Frequently studied symmetric polyhedra 
have names like ``Travelling Salesman'', ``Assignment'',
``Matching'' and ``Cut''. For these and further examples
we refer to~\cite{schrijver-2003} and the
numerous references therein. 
Over the years a rich combinatorial and geometric theory 
of polyhedra has been developed (see~\cite{ziegler-1997}, \cite{gruenbaum-2003}).
Symmetry itself is clearly a central topic in mathematics,
and through the spread of computer algebra systems like \cite{gap} and
\cite{magma}, sophisticated tools from Computational Group Theory are widely used today
(see~\cite{heo-2005}). 
Nevertheless, although many polyhedral problems 
are modeled with a high degree of symmetry, standard computational techniques 
for their solution do not take advantage of them.
Even worse, often the used methods 
are known to work notoriously poorly on symmetric problems.

\medskip

In this short survey we describe three main areas of polyhedral computations,
in which the rich geometric structure of symmetric polyhedra 
can potentially be used for improved algorithms: 

\begin{itemize}
\setlength{\itemsep}{0.0cm} 
\item[I:] Polyhedral representation conversion using symmetry
\item[II:] Symmetric integer linear programming
\item[III:] Counting lattice points and exact volumes of symmetric polyhedra
\end{itemize}

There are multiple strong dependencies among
the three topics and each one has its theoretical and algorithmic
challenges as well as important applications.
Before we take a closeup view on the three topics 
we give a brief introduction to the different types of polyhedral
symmetries and how these can be determined and worked with.

\section*{What are polyhedral symmetries?}

The symmetries of a polyhedron can be of a purely combinatorial nature
or they can also have a geometric manifestation
as {\em affine symmetries}, that is, as affine maps of $\R^n$
preserving the polyhedron. 
Among these symmetries are the ``more visually accessible'' {\em isometries} which 
are composed of translations, rotations and reflections.
All of the symmetries of the $n$-cube for example 
are part of its isometry group. 
There exists a representation as a {\em linear group} in $\GL_n (\R)$ 
and as a {\em finite orthogonal group} of isometries. 
However, if
we perturb the defining inequalities a bit, all of these affine symmetries may be lost,
while the new polyhedron is still {\em combinatorially equivalent} to
a cube, sharing all of its {\em combinatorial symmetries}. 
These are defined as automorphisms of the {\em polyhedral face lattice}
which encodes the combinatorial structure of a polyhedron.
We refer to our survey~\cite{bdprs-2012} for further
reading on these different types of polyhedral symmetry groups.
The study of combinatorial lattices and their automorphisms
is itself an active research area (see \cite{ms-2002}).
The same is true for the study of possible isometry groups, respectively of 
{\em finite orthogonal groups} in $O_n(\R)$.
Their classification becomes in a way impractical for $n\geq 5$ (see~\cite{mathoverflow37136}),
despite the classification of finite simple groups (see~\cite{atlas}).
Even less is known about symmetry groups of polyhedra (see~\cite{robertson-1984}).
Here, an ``implication phenomenon'' occurs,
which has not much been studied so far.
For instance, if a $4$-gon has an element of order~$4$ among its
affine symmetries, the $4$-gon has to be the affine image of a
square ($2$-cube), with an affine symmetry group of order~$8$.
These kind of implications clearly can potentially be exploited
algorithmically,
for example when detecting polyhedral symmetries.

\medskip

It is important to note that the same abstract group can have different 
affine representations.  
We think that a key ingredient for future algorithmic improvements 
will be the use of geometric information coming with
the affine representations of polyhedral symmetry groups. 
By a basic result in
representation theory there is an invariant affine
subspace~$\mathcal{I}$ coming with each affine symmetry group.
The polyhedron splits nicely into an invariant part $P\cap \mathcal{I}$ and
symmetric slices orthogonal to it. These lie in fibers (pre-images) of the orthogonal
projection onto~$\mathcal{I}$. In a way, all of the symmetry
is within these fibers.

\medskip

Given a polyhedron with a group of symmetries, we say
two vertices (or inequalities) are {\em equivalent}, if there exists
a group element that maps one to the other.
The set of vertices (and the set of {\em facets} / defining inequalities)
splits into a number of {\em orbits} (disjoint sets of equivalent elements).
For example, the $n$-cube has only one orbit of vertices
and one orbit of facets. 
The same is true for all Platonic polyhedra and their higher
dimensional analogues. In contrast,
the Archimedean polyhedra like
the soccer ball (truncated icosahedron) 
have more than one orbit of facets,
but only one orbit of vertices.
In all of these examples, their combinatorial symmetry group 
is equal to the group of affine symmetries.
Its invariant affine subspace is a single point, the
barycenter of the vertices.

In general, for a polyhedron~$P$ with a group of affine symmetries,
the vertices of the polyhedron split into orbits $O_1,\ldots, O_l$
and the invariant part 
$P \cap {\mathcal I}$
is equal to the convex hull
$\conv \{ b_1 , \ldots , b_l \}$,
with $b_i= (\sum_{x\in O_i} x) / |O_i|$
being the barycenter of orbit $O_i$.
This is due to two facts: 
The barycenter map, taking a point to 
the barycenter of its orbit, is an affine map.
And second, the affine image of a convex hull of given points is equal 
to the convex hull of their affine images.

Thus working with the lower dimensional polyhedron 
$P \cap {\mathcal I}$ and its vertices gives us access to 
vertices of~$P$. Orbits of integral points in~$P$ have 
barycenters at specific locations in $P \cap {\mathcal I}$. 
For instance, if the group acts transitive on the coordinates of~$\R^n$
then orbits have barycenters at integral multiples of
$\left(\frac{1}{n},\ldots,\frac{1}{n}\right)$. For more general coordinate permutations
the barycenters form a scaled copy of a standard lattice
(see~\cite{hrs-2012}).

\section*{How to determine and work with symmetries?}

If the symmetries of the polyhedron are not known,
the first difficulty is their determination and how to represent them.
In general we like to work with
as many symmetries as possible. However, the combinatorial 
symmetries can usually not be found without having full knowledge
about the vertex-facet incidences of the polyhedron (see~\cite{ks-2003}).
In contrast, the group of affine symmetries can be determined
from the vertices or defining inequalities alone, by finding the
automorphism group of an edge colored graph.
If $P$ is given as the convex hull of its vertices $x_1,\ldots , x_k$, for instance,
then the affine symmetry group 
can be obtained from the automorphism group
of the complete graph with $k$ vertices and edge labels 
$x_i^t Q^{-1} x_j$, where $Q =\sum_{i=1}^{k} x_i x_i^t$.
For details and a proof we refer to~\cite{bds-2009}.
For further methods to compute polyhedral symmetry groups we refer to~\cite{bdprs-2012}.
Automorphism groups of graphs can be computed with software like 
\cite{bliss} or \cite{nauty}.
Given a polyhedral description, the affine symmetries can 
conveniently be obtained directly with our software~\cite{sympol},
which by now can also be used through~\cite{polymake}.
For instance, given a polyhedron with its description contained in {\tt input-file}, simply call:
\begin{quote}
{\tt sympol --automorphisms-only  input-file}
\end{quote}

\medskip

If the symmetry group of a polyhedron (or parts of it) are given as a {\em permutation group},
we can use sophisticated tools from {\em Computational Group Theory}.
Each element of the group is then viewed as a {\em permutation} of the
index set $\{1,\ldots,m\}$ of the input, for instance of $m$ defining inequalities.
In practice, it is necessary to work with a small set of {\em group generators}
if the group is large, and there are advanced heuristics to obtain such
sets.
Each face (and in particular each vertex) of a polyhedron is determined by a number of inequalities
that are satisfied with equality; it can therefore be represented
by a subset of $\{1,\ldots,m\}$.
Given generators of a large permutation group
and two subsets that represent faces, 
a typical computational bottleneck is to decide 
if both are in the same orbit. 
The fundamental data structures used for this in practice
are {\em bases and strong generating sets}
(BSGS, see \cite{seress-2003}, \cite{heo-2005}).
Based on them, backtrack searches can be used to perform essential tasks,
such as deciding on (non-)equivalence, obtaining stabilizers or
fusing and splitting orbits.
An elaborate version is the {\em partition backtrack} introduced
by Leon  \cite{leon-1991}.
These backtracking methods work quite well in practice, although
from a complexity point of view the mentioned problems are thought to be difficult
(see \cite{luks-1993}).
Although computer algebra systems like~\cite{gap}
and~\cite{magma}
provide functions to work with permutation groups,
for performance reasons it is often
desirable to use problem specific code (see for example~\cite{ko-2006}).
Nevertheless, all of these approaches, 
including~\cite{gap} and~\cite{magma},
rely on efficient implementations of some partition backtrack. 
Therefore we have created a flexible C++-implementation \cite{permlib}
of Leon's partition backtrack (see \cite{rehn-2010a})
that can serve as a basis for the development of algorithms 
which combine tools from Computational Group Theory and Polyhedral Combinatorics.
By now, \cite{permlib} has successfully been integrated into 
current versions of \cite{sympol}, \cite{polymake} and \cite{scip}.

\section*{I. Representation Conversion}

By a fundamental theorem in polyhedral combinatorics, 
the Farkas-Minkowski-Weyl theorem,
every polyhedron has a second representation as the convex hull 
of finitely many {\em vertices}
(extreme points) and, in the unbounded case, some {\em rays} 
(see~\cite{ziegler-1997}).
Converting representations from inequalities to vertices (and rays), or 
vice versa, is a frequent task known as {\em representation conversion problem}
(or {\em convex hull problem}).
The importance of these conversions is due to the fact that some
problems, like the maximization of a nonlinear convex function, 
are easy to solve in one presentation, but not in the other.
Often, vertices represent objects that one would like to 
classify. These objects can be quite diverse, for instance 
perfect quadratic forms (see~\cite{dsv-2007})
or the elementary flux modes in biochemical reaction systems (see~\cite{sh-1994}).
Representation conversions are also often used to analyze
polyhedra in Combinatorial Optimization (see~\cite{schrijver-2003}).
So far there exists no efficient algorithm for finding all the vertices of a polyhedron.
In fact, the existence of such an algorithm appears to be unlikely, 
as it is NP-complete for polyhedra that are possibly unbounded (see~\cite{kbbeg-2008}).
Nevertheless, several algorithms and implementations are
widely used in practice
(see for example~\cite{cdd} and~\cite{lrs}
which are also available through~\cite{polymake}).

\medskip

Quite often one is only interested in one representative
for each orbit of vertices (or inequalities) in a representation conversion.
For example, when maximizing a nonlinear
convex function on a polyhedron, or when vertices and inequalities 
in one orbit correspond to equivalent objects of some sort.
{\em Representation conversion up to symmetries}
has been considered in different contexts, and depending on the problem,
different techniques have been successful. 
The most successful approaches currently known are the 
{\em Incidence Decomposition Method}
and the {\em Adjacency Decomposition Method}
(see for instance \cite{cr-1996}, \cite{dfps-2001}, \cite{dfmv-2003},
\cite{dsv-2007}, \cite{dsv-2009}, \cite{dsv-2010}).
Both methods decompose the problem into a number of lower dimensional
subproblems. They can be used recursively and can be parallelized
(see~\cite{cr-2002}, \cite{di-2007}).
Loosely speaking, the Incidence Decomposition Method fixes an orbit of
the input, whereas the Adjacency Decomposition Method fixes an orbit of the
output and then lists all ``neighboring'' orbits.
For details we refer to our survey~\cite{bds-2009}.
We note that it is a priori not clear which method works best.
We think best results can be achieved by a combination of
different algorithms.
All methods known so far do not use geometric insights and still rely on
subproblem conversions that do not exploit available symmetry.

\medskip

{\bf Our software~\cite{sympol}} 
and the experimental \cite{gap} package~\cite{polyhedral}
provide implementations of decomposition methods. These preliminary tools have
already successfully been used in our own work
(\cite{dsv-2007}, \cite{bbcgks-2009}, \cite{dsv-2009}, \cite{dsv-2010}, \cite{dse-2011}), 
but also by others: For instance, Kumar~\cite{kumar-2011} obtains
a classification of elliptic fibrations that was previously impossible.
Jacques Martinet writes in \cite{martinet-2003} about the result in \cite{dsv-2007}:
``It seems plainly impossible to classify $8$-dimensional perfect lattices.''
\cite{sympol} can also be used to verify cumbersome
calculations in proofs, like
the edge-graph diameter analysis of the recently discovered, celebrated 
counterexamples to the Hirsch conjecture~(see~\cite{santos-2012}, \cite{rs-2010}). 
For example, with the call
\begin{quote}
{\tt sympol --idm-adm-level 0 1  --adjacencies  input-file}
\end{quote}
where {\tt input-file} contains the 48 vertices of the $5$-dimensional {\em Santos prismatoid}
(see Table~1 in \cite{santos-2012}),
{\tt SymPol} returns a text file with a description of the adjacency graph of facets 
up to symmetry. Using a visualization tool like \cite{graphviz}, the
produced textfile, say {\tt adjacencies.dot}, 
can then easily be turned into an image with a command like
\begin{quote}
{\tt neato -Tpng -o adjacencies.png adjacencies.dot}
\end{quote}
From the obtained image (see figure) 
it is easily verified that the shortest path from facet~1 to facet~12
is of length $6$, which is the key calculation in the proof of~\cite{santos-2012}.

\medskip

\centerline{\includegraphics[width=5cm]{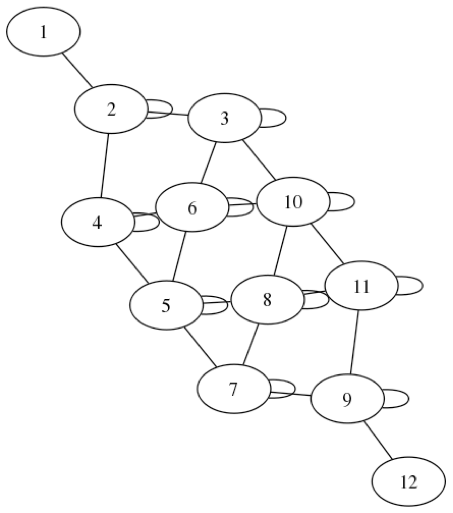}}

\medskip

Let us make a remark on the increasing importance of mathematical software in general:
As sophisticated computational tools become an increasingly important basis 
for high-level mathematical research, their creation 
also becomes an increasingly important service to the mathematical community. 
More and more mathematicians use computers in the creative process 
and to verify standard parts of difficult larger proofs (see
e.g.~\cite{hales-2005},~\cite{ck-2009},~\cite{hales-2011}).
Timothy Gowers~\cite{gowers-2000} even guesses that at the end of the
21st century, 
computers will be better than humans in proving theorems.
Although we would not go as far as Gowers, 
we are convinced that in the future, parts of proofs will routinely be
performed by computers.
With a symbiosis of human and computer reasoning we will see 
substantial advances in mathematical problems. 
In this way reliable mathematical software becomes an increasingly 
important part of mathematics itself.

\subsection*{Challenges} 
One of the most challenging polyhedral conversion problems 
arises in conjunction with lattice sphere packing problem, a classical problem in
the Geometry of Numbers. 
Since its solution up to dimension~$8$ almost $80$ years
ago, it is still open in dimension $n\geq 9$, 
with the exception of dimension $24$ (see~\cite{martinet-2003}, \cite{ck-2009}).
One way to approach this problem is
via finding the vertices of a locally polyhedral object
known as  {\em Ryshkov polyhedron} (see \cite{schuermann-2009a} for details).
The currently open $9$-dimensional case leads to a challenging 
representation conversion problem of a $45$-dimensional Ryshkov polyhedron.
Main difficulties here come from faces that
carry the symmetries of the exceptional 
Weyl group~$\mathsf{E}_8$.
We think that this problem is a 
particular nice test case, as all finite rational matrix groups appear 
as stabilizers of faces in the Ryshkov polyhedron.
So in a way, this challenging representation conversion problem
gives a {\em universal test case} for any future algorithmic advances.

\section*{II. Integer Linear Programming}

{\em Linear programming} is the task of maximizing (or minimizing)
a linear function on a polyhedron given by linear inequalities. 
It serves as a fundamental basis for theory and computations in Mathematical
Optimization~(see~\cite{todd-2002}).
In {\em Integer linear programming} vectors to be optimized
are restricted to integers (number of goods, etc.)
or even to $0/1$ entries (encoding a simple yes-no-decision).
Integer linear programming is widely used in practical applications.
In fact,
``the vast majority of applications found in operations research and industrial engineering involve 
the use of discrete variables in problem formulation'' (from a
book review of \cite{wolsey-1998}).
In many of these problems the involved polyhedra have  
symmetries (see \cite{schrijver-2003}).
From a complexity point of view, integer programming is
NP-complete (see \cite{karp-1972}),
whereas linear programming can be solved in polynomial time.  
For fixed dimension,
polynomial time algorithms are known for integer linear programming
(see~\cite{lenstra-1983}).

\medskip

A linear programming problem $\max \; c^t x  \; \mbox{with} \; x \in P$ 
is {\em invariant with respect to a linear symmetry group} $\Gamma\leq \GL_n(\R)$, 
if the polyhedron $P$ and the utility vector $c\in\R^n$ are preserved by it,
that is, if $\Gamma P = P$ and  $\Gamma c = c$.
Any solution to the linear program, its orbit and the barycenter of its orbit 
lie in the same hyperplane orthogonal to the utility vector and therefore have the same 
utility value.
Due to the convexity of polyhedra, the
barycenter is also a feasible solution. As it lies in the 
{\em invariant linear subspace} $\mathcal{I}$ of ${\Gamma}$,
the linear programming problem always has a solution 
attained within~$\mathcal{I}$.
Thus it is possible to solve the lower-dimensional 
linear program $\max \; c^t x  \; \mbox{with} \; x \in P\cap
\mathcal{I}$.
Such symmetry reductions are often referred to as ``dimension reduction''
or ``variable reduction''.
The symmetries of an integer linear program are more
restrictive, as also $\Z^n$ has to be left invariant by the group~$\Gamma$.

\medskip

Exploiting symmetries in integer programming is much more difficult
than in linear programming.
In fact, symmetries are rather problematic, as
standard methods like {\em branch-and-bound} 
or {\em branch-and-cut} (see \cite{schrijver-1986})
have to solve many equivalent sub-problems in such cases.
In contrast to linear programming, it is not possible to simply
consider the intersection with the invariant affine subspace, 
as integral solutions can lie outside.
Nevertheless, in recent years it has been shown
that it is possible to exploit symmetries in integer programming;
see for example
\cite{margot-2003}, \cite{friedman-2007}, \cite{bm-2008},
\cite{kp-2008}, \cite{olrs-2011}, \cite{lmt-2009}.
These specific methods fall into two main classes:
They either modify the standard branching approach,
using isomorphism tests or isomorphism free generation to
avoid solving equivalent subproblems; or they use techniques
to cut down the original symmetric problem to a less symmetric one,
which contains at least one element of each orbit of solutions.
For further reading we refer to the excellent survey~\cite{margot-2009}.


As many real world applications can be modeled as
(mixed) integer programming problems, a variety of 
professional software packages are available. 
Two of the leading ones,
\cite{cplex} and \cite{gurobi}, have by now included some 
techniques to avoid or use symmetry.
Unfortunately it is publicly not known what exactly is done.

\medskip

None of the known methods uses the rich geometric 
properties of the involved symmetric polyhedra.
Using the fact that solutions are ``near'' the linear invariant
subspace, it is possible to do better.
For the special case of
a one dimensional invariant subspace,
with the full symmetric group $S_n$ 
acting transitively on the coordinates in $\R^n$,
this is shown in~\cite{bhj-2011}. 
We have highly promising
results with a generalization to arbitrary symmetries in~\cite{hrs-2012}:
In particular for direct products of symmetric groups,
we not only beat state-of-the-art professional solvers, 
but even solve a challenging, previously
unsolved benchmark problem from~\cite{miplib}
(instance {\tt toll-like}).

The main ingredient is the observation that any 
feasible integer linear programming problem with a non-trivial affine symmetry group
contains certain {\em core-sets} of integral vectors that can be used as a 
kind of test-set, that is, if non of the points from the core-set is contained 
in the feasible region, no integral point is. 
Assume $\Gamma\leq \GL_n(\R)$ is a linear representation of a given symmetry group
and ${\mathcal{I}}$ denotes its invariant linear subspace containing the utility vector~$c$.
Then we say an integral point~$z$ in a  {\em fiber} (pre-image) of the orthogonal
projection onto ${\mathcal{I}}$
is in the {\em core-set of the fiber} if the convex hull 
of its orbit $\Gamma z$ does not contain any integral points aside
those of the orbit itself. 
Then, by definition, representatives of each orbit $\Gamma z$ in the
core-set can be used as a {\em test-set for feasibility of a fiber}.

Using the Flatness Theorem (from the Geometry of Numbers), 
it can be shown that core-sets are finite for irreducible groups.
For direct products of symmetric groups 
acting on some of the coordinates, 
the test-set containing only representatives of orbits 
even reduces to a single point.
Besides that not much is known so far about core-sets.
Nevertheless, we think that they will serve as a powerful tool in the design 
of new algorithms for symmetric integer linear programs.


%

\subsection*{Challenges} 

Challenging examples of symmetric integer linear programs can be found
in benchmark libraries like~\cite{miplib}. These problems come from diverse 
contexts and have not been chosen to be particular symmetric.
Nevertheless many symmetries can be found and exploiting them algorithmically,
beating state-of-the-art commercial solvers, remains a challenging test case
for future advances.

Some particular symmetric integer linear programming problems 
coming from difficult combinatorial problems in mathematics 
have been collected (and worked on) by Francois Margot at \cite{symlplib}
(see for instance also \cite{margot-2003} and \cite{bm-2008}).
As these problems have been intensively worked on, 
improving on the currently known results is
certainly a hard problem. So this gives a very good benchmark for future 
improvements as well.

\section*{III. Lattice Point Counting and Exact Volumes}

Often it is desirable to know how many integral solutions there are to
a system of linear inequalities. Such problems occur frequently in
Combinatorics (see~\cite{stanley-1997}) 
but also in disciplines such as 
representation theory (Kostka and Littlewood-Richardson
coefficients, see \cite{bz-2001} and \cite{kt-1999}),
in statistics (contingency tables, see \cite{ds-1998}),
in voting theory (see \cite{wp-2007} and \cite{gl-2011}), 
and even in compiler optimization (see~\cite{graphite}). 
We refer to \cite{deloera-2005} for an overview.
Counting lattice points is moreover 
intimately related to integer linear programming (see~\cite{lasserre-2009}).

\medskip

By a breakthrough result of Barvinok~\cite{barvinok-1994} in the
1990s, counting lattice points inside a rational polyhedron can be
done in polynomial time for fixed dimension. His ideas
are based on evaluating ``short rational generating functions''
and on constructing unimodular triangulations.
His algorithm has been implemented in \cite{latte} and \cite{barvinok}.
The same applies to a slightly more general setting, in which one
considers a one-parameter family of dilations $\lambda P$, with $P$ a
rational polyhedron and $\lambda$ an integer. By a theory initiated by Ehrhart in 
the 1960s (see \cite{ehrhart-1967}, \cite{br-2007}), it is possible to obtain the
number of integral points in the dilate $\lambda P$ 
by a quasi-polynomial in $\lambda$, with its degree equal to the
dimension of~$P$.
A quasi-polynomial~$p$ is determined by a finite number of polynomials 
$p_i$, $i=0,\ldots, k$,
via the setting $p(\lambda)=p_i(\lambda)$ for all~$\lambda$ congruent to~$i$ mod $k$.
In case of $P$ being integral, the {\em Ehrhart quasi-polynomial} simply 
is a polynomial in~$\lambda$.
In general, 
the quasi-polynomial can also be computed in polynomial time by Barvinok methods 
(see~\cite{barvinok-2008}).
Often, the main interest is only in the leading coefficient of the Ehrhart
quasi-polynomial, which is the volume of~$P$.
Computing the volume itself is already a $\#P$-hard problem (see
\cite{df-1988}, \cite{bw-1991}).

\medskip

Despite the fact that many 
counting problems have plenty of symmetries,
they have not been exploited systematically so far.
In other words,
exploiting symmetry in lattice point counting, or more generally in
Ehrhart theory, is a vastly open subject.
For volume computations the situation seems a bit better. For very special 
volume computations symmetry can be exploited (see \cite{dsv-2009}).
However, there is still a huge potential for improved methods.
Many of the difficulties
originate from the fact that the ``Barvinok methods'' used to solve
them rely on unimodular triangulations of polyhedral cones that 
usually do not inherit the symmetry of the polyhedron. 
New roads will have to be taken here.

\medskip

In~\cite{schuermann-2011} it is shown that 
it is possible to exploit symmetries 
by using a decomposition into symmetric slices, together with a 
{\em weighted Ehrhart theory}.
The theoretical background and
first implementations for such a theory have just recently been developed
(see~\cite{bbdmv-2011b}, \cite{bs-2012}).
The polyhedral decomposition used in~\cite{schuermann-2011}  is rather special: 
There is a linear invariant part and symmetric slices orthogonal to it, which are cross-products of
regular simplices (simplotopes). 
A generalization to other decompositions is easily obtained, whenever there is 
a decomposition into an invariant part and slices orthogonal to it
for which the Ehrhart quasi-polynomial is known. 
Note that the decomposition can easily be obtained in a automated way,
as the invariant part is the intersection of the given polyhedron with
the affine space fixed by its symmetry group.

For exploiting symmetry in corresponding volume computations, the 
integration of polynomials over a polyhedron is used.
Using Brion-Lawrence-Varchenko theory, this can efficiently be done
by integrating sums of powers of linear forms (see~\cite{bbdmv-2011a}.
The new decomposition approach of~\cite{schuermann-2011}
also allows to obtain
exact volumes that have not been computable before.
This is demonstrated on three well studied examples from Social Choice theory,
which give the exact likelihood of certain election outcomes with four candidates 
that were previously known for three candidate elections only
(see~\cite{gl-2011}).

\subsection*{Challenges} 

In Social Choice theory we face a large amount of challenging problems
related to probability calculations of voting situations with four and more candidates.
The only known results in the context of the ``polyhedral model'' (IAC
hypothesis) appear to be those in~\cite{schuermann-2011}, which are obtained by
exploiting polyhedral symmetry as described above. 

A challenging benchmark volume computation that several researchers previously have looked
at is the volume of the Birkhoff polytope $B_n$
(also known as {\em perfect matching polytope} of the complete bipartite graph $K_{n,n}$).
The current known record is the volume of~$B_{10}$ 
due to Beck and Pixton \cite{bp-2003}, using a 
complex-analytic way to compute the Ehrhart polynomial. 
The computation of the
volume of~$B_{11}$ would certainly be quite a computational achievement.

\section*{Conclusions}

We expect that symmetry exploiting techniques for polyhedral computations 
can be vastly improved by using geometric properties that come with affine symmetries of 
polyhedra. Concentrating on improvements in polyhedral computations
with affine symmetries is practically no restriction:
If a polyhedron is given, either by linear inequalities or vertices and rays, 
the affine symmetries of the (potentially larger)
combinatorial symmetry group 
are practically the only ones we can compute.

For polyhedral representation conversions we see potential in
enhancing decomposition methods through the use of geometric
information like fundamental domains, classical invariant theory 
and symmetric polyhedral decompositions.
For integer linear programming we expect that a new class of algorithms based on
the concept of core points will help to exploit symmetry on difficult symmetric
integer linear programming problems.
For exact volume computations and counting of lattice points, there is still 
a lot of potential for new ideas using symmetry.
Overall, we think 
symmetry should be exploitable whenever it is available.
For this goal to be reached there seem still quite some efforts necessary though.

\newcommand{\etalchar}[1]{$^{#1}$}
\providecommand{\bysame}{\leavevmode\hbox to3em{\hrulefill}\thinspace}
\providecommand{\href}[2]{#2}

\renewcommand\refname{}  

\end{document}